\documentclass{article}
\usepackage{amssymb}
\usepackage{amsmath}
\usepackage{amsfonts}

\newtheorem{theorem}{Theorem}

\newtheorem{lemma}{Lemma}

\newtheorem{remark}{Remark}

\newdimen\dummy
\dummy=\oddsidemargin
\begin{document}

\title{\textbf{On the degree of strong approximation of almost periodic
functions in the Stepanov sense}}
\author{\textbf{W\l odzimierz \L enski and Bogdan Szal} \\
University of Zielona G\'{o}ra\\
Faculty of Mathematics, Computer Science and Econometrics\\
65-516 Zielona G\'{o}ra, ul. Szafrana 4a, Poland\\
W.Lenski@wmie.uz.zgora.pl, B.Szal @wmie.uz.zgora.pl }
\date{ \ \ }
\maketitle

\begin{abstract}
Considering the class of almost periodic functions in the Stepanov sense we
extend and generalize the results of the first author \cite{4}. as well as
the results of L. Leindler \cite{3} and P. Chandra \cite{1,2} .


\ \ \ \ \ \ \ \ \ \ \ \ \ \ \ \ \ \ \ \ 

\ \ \ \ \ \ \ \ \ \ \ \ \ \ \ \ \ \ 

\end{abstract}

\section{Introduction}

Let $S^{p}\;\left( 1<p\leq \infty \right) $ be the class of all almost
periodic functions in the Stepanov sense $\left( 1<p<\infty \right) $ or
uniformly almost periodic $\left( p=\infty \right) $ with the norm 
\begin{equation*}
\Vert f\Vert _{S^{p}}:=\left\{ 
\begin{array}{c}
\sup\limits_{u}\left\{ \frac{1}{\pi }\int_{u}^{u+\pi }\mid f(t)\mid
^{p}dt\right\} ^{1/p}\text{ \ \ when \ \ }1<p<\infty , \\ 
\sup\limits_{u}\mid f(u)\mid \text{ \ \ when \ \ }p=\infty .%
\end{array}%
\right.
\end{equation*}%
Suppose that the Fourier series of $f\in S^{p}$ has the form 
\begin{equation*}
Sf\left( x\right) =\sum_{\nu =-\infty }^{\infty }A_{\nu }\left( f\right)
e^{i\lambda _{\nu }x},\text{ \ \ where \ }A_{\nu }\left( f\right)
=\lim_{L\rightarrow \infty }\frac{1}{L}\int_{0}^{L}f(t)e^{-i\lambda _{\nu
}t}dt,
\end{equation*}%
with the partial sums\ 
\begin{equation*}
S_{\gamma _{k}}f\left( x\right) =\sum_{\left\vert \lambda _{\nu }\right\vert
\leq \gamma _{k}}A_{\nu }\left( f\right) e^{i\lambda _{\nu }x}
\end{equation*}%
and that $0=\lambda _{0}<\lambda _{\nu }<\lambda _{\nu +1}$ if $\nu \in 
\mathbb{N}
=\left\{ 1,2,3...\right\} ,$ $\underset{v\rightarrow \infty }{\lim }\lambda
_{\nu }=\infty ,$ $\lambda _{-\nu }=-\lambda _{\nu ,}$\ $\left\vert A_{\nu
}\right\vert +\left\vert A_{-\nu }\right\vert >0.$ Let\ $\Omega _{\alpha ,p}$%
, with some fixed positive $\alpha $ , be the set of functions of class $S^{p%
\text{ }}$ whose Fourier exponents satisfy the condition 
\begin{equation*}
\lambda _{\nu +1}-\lambda _{\nu }\geq \alpha \text{ \ \ }\left( \nu \in 
\mathbb{N}
\right) .
\end{equation*}%
In case $f\in \Omega _{\alpha ,p}$\ 
\begin{equation*}
S_{\lambda _{k}}f\left( x\right) =\int_{0}^{\infty }\left\{ f\left(
x+t\right) +f\left( x-t\right) \right\} \Psi _{\lambda _{k},\lambda
_{k}+\alpha }\left( t\right) dt,
\end{equation*}%
where 
\begin{equation*}
\Psi _{\lambda ,\eta }\left( t\right) =\frac{2\sin \frac{\left( \eta
-\lambda \right) t}{2}\sin \frac{\left( \eta +\lambda \right) t}{2}}{\pi
\left( \eta -\lambda \right) t^{2}}\text{ \ \ }\left( 0<\lambda <\eta ,\text{
\ }\left\vert t\right\vert >0\right) .
\end{equation*}

Let $A:=\left( a_{nk}\right) $ $\left( k,n=0,1,...\right) $ be a lower
triangular infinite matrix of real numbers satisfying the following
condition:

\begin{equation}
\text{ }a_{nk}\geq 0\text{ }\left( k,n=0,1,...\right) ,\text{ }a_{nk}=0\text{
}\left( k>n\right) \text{ and }\sum\limits_{k=0}^{n}a_{nk}=1.  \label{1}
\end{equation}

Let us consider the strong mean 
\begin{equation}
H_{n,A,\gamma }^{q}f\left( x\right) =\left\{ \sum_{k=0}^{n}a_{n,k}\left\vert
S_{\gamma _{k}}f\left( x\right) -f\left( x\right) \right\vert ^{q}\right\}
^{1/q}\text{ \ \ \ }\left( q>0\right) \text{.}  \label{S1}
\end{equation}%
As measures of approximation by the quantity (\ref{S1}), we use the best
approximation of $f$ by entire functions $g_{\sigma }$ of exponential type $%
\sigma $\ bounded on the real axis, shortly $g_{\sigma }\in B_{\sigma }$ and
the moduli of continuity of $\ f$ defined by the formulas%
\begin{equation*}
E_{\sigma }(f)_{S^{p}}=\inf_{g_{\sigma }}\left\Vert f-g_{\sigma }\right\Vert
_{S^{p}},
\end{equation*}%
\begin{equation*}
\omega f\left( \delta \right) _{S^{p}}=\sup_{\left\vert t\right\vert \leq
\delta }\left\Vert f\left( \cdot +t\right) -f\left( \cdot \right)
\right\Vert _{S^{p}\text{ }}
\end{equation*}%
and%
\begin{equation*}
w_{x}f(\delta )_{p}:=\left\{ \frac{1}{\delta }\int_{0}^{\delta }\left\vert
\varphi _{x}\left( t\right) \right\vert ^{p}dt\right\} ^{1/p}\text{ with }%
1<p<\infty ,
\end{equation*}%
where $\varphi _{x}\left( t\right) :=f\left( x+t\right) +f\left( x-t\right)
-2f\left( x\right) $, respectively.

A sequence $c:=\left( c_{n}\right) $ of nonnegative numbers tending to zero
is called the Rest Bounded Variation Sequence, or briefly $c\in RBVS$, if it
has the property%
\begin{equation}
\sum\limits_{k=m}^{\infty }\left\vert c_{n}-c_{n+1}\right\vert \leq K\left(
c\right) c_{m}  \label{2}
\end{equation}%
for all natural numbers $m$, where $K\left( c\right) $ is a constant
depending only on $c$.

A sequence $c:=\left( c_{n}\right) $ of nonnegative numbers will be called
the Head Bounded Variation Sequence, or briefly $c\in HBVS$, if it has the
property%
\begin{equation}
\sum\limits_{k=0}^{m-1}\left\vert c_{n}-c_{n+1}\right\vert \leq K\left(
c\right) c_{m}  \label{3}
\end{equation}%
for all natural numbers $m$, or only for all $m\leq N$ if the sequence $c$
has only finite nonzero terms and the last nonzero terms is $c_{N}$.

Therefore we assume that the sequence $\left( K\left( \alpha _{n}\right)
\right) _{n=0}^{\infty }$ is bounded, that is, that there exists a constant $%
K$ such that%
\begin{equation*}
0\leq K\left( \alpha _{n}\right) \leq K
\end{equation*}%
holds for all $n$, where $K\left( \alpha _{n}\right) $ denote the sequence
of constants appearing in the inequalities (\ref{2}) or (\ref{3}) for the
sequence $\alpha _{n}:=\left( a_{nk}\right) _{k=0}^{\infty }$.Now we can
give the conditions to be used later on. We assume that for all $n$ and $%
0\leq m\leq n$%
\begin{equation}
\sum\limits_{k=m}^{\infty }\left\vert a_{nk}-a_{nk+1}\right\vert \leq
Ka_{nm}  \label{4}
\end{equation}%
and%
\begin{equation}
\sum\limits_{k=0}^{m-1}\left\vert a_{nk}-a_{nk+1}\right\vert \leq Ka_{nm}
\label{5}
\end{equation}%
hold if $\alpha _{n}:=\left( a_{nk}\right) _{k=0}^{\infty }$ belongs to $%
RBVS $ or $HBVS$, respectively.

The $C$-norm of the deviation $\left\vert \sum_{k=0}^{n}a_{n,k}\left[
S_{k}f\left( x\right) -f\left( x\right) \right] \right\vert ,$ with the
partial sums $S_{k}f$ of classical trigonometric Fourier series, was
estimated by P. Chandra \cite{1} \cite{2} for monotonic sequences $\left(
a_{nk}\right) $ and by L. Leindler \cite{3} for the sequences of bounded
variation. These results were generalized by W. \L enski \cite{4} who
considered the strong means $H_{n,A}^{q\text{ }},$ also in classical case,
and the functions belonging to the $L^{p}$. In present paper we shall
considered the almost periodic functions from the Stepanov class giving
similarly estimations for the strong means $H_{n,A}^{q\text{ }}$ in
individual points and in norms.

We shall write $I_{1}\ll I_{2}$ if there exists a positive constant $C$ such
that $I_{1}\leq CI_{2}$.

\section{Main results \ \ \ \ }

Let us consider a function $w_{x}$ of modulus of continuity type on the
interval $[0,+\infty ),$ i.e. a nondecreasing continuous function having the
following properties:\ $w_{x}\left( 0\right) =0,$ $w_{x}\left( \delta
_{1}+\delta _{2}\right) \leq w_{x}\left( \delta _{1}\right) +w_{x}\left(
\delta _{2}\right) $ for any $\delta _{1},\delta _{2}\geq 0$\ with $x$ such
that the set%
\begin{eqnarray*}
\Omega _{\alpha ,p}\left( w_{x}\right) &=&\left\{ f\in \Omega _{\alpha ,p}: 
\left[ \frac{1}{\delta }\int_{0}^{\delta }\left\vert \varphi _{x}\left(
t\right) -\varphi _{x}\left( t\pm \gamma \right) \right\vert ^{p}dt\right]
^{1/p}\ll w_{x}\left( \gamma \right) \right. \\
&&\left. \text{\ and \ }w_{x}f\left( \delta \right) _{p}\ll w_{x}\left(
\delta \right) \text{ \ , \ where \ }\gamma ,\delta >0\right\}
\end{eqnarray*}%
is nonempty. It is clear that $\Omega _{\alpha ,p}\left( w_{x}\right)
\subseteq \Omega _{\alpha ,p^{\prime }}\left( w_{x}\right) ,$ for $p^{\prime
}\leq p<\infty .$

Our main results are the following:

\begin{theorem}
Let (\ref{1}) and (\ref{5}) hold. Suppose $w_{x}$ is such that 
\begin{equation}
\left\{ u^{\frac{p}{q}}\int\limits_{u}^{\pi }\frac{\left( w_{x}\left(
t\right) \right) ^{p}}{t^{1+\frac{p}{q}}}dt\right\} ^{\frac{1}{p}}=O\left(
uH_{x}\left( u\right) \right) \text{ \ \ as \ \ }u\rightarrow 0^{+},
\label{6}
\end{equation}%
where $H_{x}\left( u\right) \geq 0$, $1<p\leq q$ and 
\begin{equation}
\int\limits_{0}^{t}H_{x}\left( u\right) du=O\left( tH_{x}\left( t\right)
\right) \text{ \ \ as \ \ }t\rightarrow 0^{+}.  \label{7}
\end{equation}%
If $f\in \Omega _{\alpha ,p}\left( w_{x}\right) ,$ then%
\begin{equation}
H_{n,A,\gamma }^{q}f\left( x\right) =O\left( a_{nn}H_{x}\left( a_{nn}\right)
+\left\{ \sum_{k=0}^{n}a_{n,k}\left( E_{\alpha k/2}\left( f\right)
_{S^{p}}\right) ^{q}\right\} ^{1/q}\right) ,  \label{8}
\end{equation}%
where $q$ is such that $1<q\left( q-1\right) ^{-1}\leq p\leq q.$
\end{theorem}

\begin{theorem}
Let (\ref{1}), (\ref{4}), (\ref{6}) and (\ref{7}) hold. If $f\in \Omega
_{\alpha ,p}\left( w_{x}\right) $, then%
\begin{equation}
H_{n,A,\gamma }^{q}f\left( x\right) =O\left( a_{n0}H_{x}\left( a_{n0}\right)
+\left\{ \sum_{k=0}^{n}a_{n,k}\left( E_{\alpha k/2}\left( f\right)
_{S^{p}}\right) ^{q}\right\} ^{1/q}\right) ,  \label{10}
\end{equation}%
where $q$ is such that $1<q\left( q-1\right) ^{-1}\leq p\leq q.$
\end{theorem}

Consequently, we can immediately derive the results on norm approximation.

\begin{theorem}
Let (\ref{1}) and (\ref{5}) hold. Suppose $\omega f\left( \cdot \right) _{S^{%
\widetilde{p}}}$ is such that 
\begin{equation}
\left\{ u^{\frac{p}{q}}\int\limits_{u}^{\pi }\frac{\left( \omega f\left(
t\right) _{S^{\widetilde{p}}}\right) ^{p}}{t^{1+\frac{p}{q}}}dt\right\} ^{%
\frac{1}{p}}=O\left( uH\left( u\right) \right) \text{ \ \ as \ \ }%
u\rightarrow 0^{+}  \label{11a}
\end{equation}%
holds, with $1<p\leq q\leq \widetilde{p}$ , where additionally $H$ $\left(
\geq 0\right) $ instead of $H_{x}$ satisfies the condition (\ref{7}). If $%
f\in \Omega _{\alpha ,\widetilde{p}}$, then%
\begin{equation*}
\left\Vert H_{n,A,\gamma }^{q^{\prime }}f\left( \cdot \right) \right\Vert
_{S^{\widetilde{p}}}=O\left( a_{nn}H_{x}\left( a_{nn}\right) \right) ,
\end{equation*}%
with $q^{\prime }\in (0,q]$, where $q$ is such that $1<q\left( q-1\right)
^{-1}\leq p\leq q.$
\end{theorem}

\begin{theorem}
Let (\ref{1}) and (\ref{4}) hold. Suppose $\omega f\left( \cdot \right) _{S^{%
\widetilde{p}}}$ is such that (\ref{11a}) holds, with $1<p\leq q\leq 
\widetilde{p}$ , where additionally $H$ $\left( \geq 0\right) $ instead of $%
H_{x}$ satisfies the condition (\ref{7}). If $f\in \Omega _{\alpha ,%
\widetilde{p}}$, then%
\begin{equation*}
\left\Vert H_{n,A,\gamma }^{q^{\prime }}f\left( \cdot \right) \right\Vert
_{S^{\widetilde{p}}}=O\left( a_{n0}H_{x}\left( a_{n0}\right) \right) ,
\end{equation*}%
with $q^{\prime }\in (0,q]$, where $q$ is such that $1<q\left( q-1\right)
^{-1}\leq p\leq q.$
\end{theorem}

\begin{remark}
Analyzing\ our proofs and dividing the integral in the formula%
\begin{equation*}
\left\{ \sum_{k=0}^{n}a_{n,k}\left\vert \int_{0}^{\infty }\varphi _{x}\left(
t\right) \Psi _{k+\kappa }\left( t\right) dt\right\vert ^{q}\right\} ^{1/q}
\end{equation*}%
into parts with $\frac{\pi }{n+1}$ instead of $a_{n,n}$ or $a_{n,0}$ we can
obtain the next series of theorems analogously as in \cite{4}.
\end{remark}

\section{Lemmas}

To prove our theorems we need the following lemmas.

\begin{lemma}
\cite{4} If (\ref{6}) and (\ref{7}) hold, then%
\begin{equation}
\int\limits_{0}^{u}\frac{w_{x}f\left( t\right) }{t}dt=O\left( uH_{x}\left(
u\right) \right) \text{ \ }\left( u\rightarrow 0_{+}\right) .  \label{12}
\end{equation}
\end{lemma}

\begin{lemma}
\cite[Theorem 5.20 II, Ch. XII]{5} Suppose that $1<q\left( q-1\right)
^{-1}\leq p\leq q$ and $\xi =\frac{1}{p}+\frac{1}{q}-1$. If $\left\vert
t^{-\xi }g\left( t\right) \right\vert \in L^{p,}$ then%
\begin{equation}
\left\{ \frac{\left\vert a_{0}\left( g\right) \right\vert ^{q}}{2}%
+\sum\limits_{k=0}^{\infty }\left( \left\vert a_{k}\left( g\right)
\right\vert ^{q}+\left\vert b_{k}\left( g\right) \right\vert ^{q}\right)
\right\} ^{\frac{1}{q}}\ll \left\{ \int\limits_{-\pi }^{\pi }\left\vert
t^{-\xi }g\left( t\right) \right\vert ^{p}dt\right\} ^{\frac{1}{p}}.
\label{13}
\end{equation}
\end{lemma}

\section{Proofs of the Results}

\subsection{Proof of Theorem 1}

In the proof we will use the following function $\Phi _{x}f\left( \delta
,\nu \right) =\frac{1}{\delta }\int_{\nu }^{\nu +\delta }\varphi _{x}\left(
u\right) du,$ with $\delta =\delta _{n}=\frac{\pi }{n+1}$ and its estimate
from \cite[Lemma 1, p.218]{WL}%
\begin{equation}
\left\vert \Phi _{x}f\left( \xi _{1},\xi _{2}\right) \right\vert \leq
w_{x}\left( \xi _{1}\right) +w_{x}\left( \xi _{2}\right)  \label{W}
\end{equation}%
for $f\in \Omega _{\alpha ,p}\left( w_{x}\right) $ and any $\xi _{1},\xi
_{2}>0$.

Since, for $n=0$ our estimate is evident we consider $n>0$, only.

Denote by $S_{k}^{\ast }f$ the sums of the form 
\begin{equation*}
S_{\frac{\alpha k}{2}}f\left( x\right) =\sum_{\left\vert \lambda _{\nu
}\right\vert \leq \frac{\alpha k}{2}}A_{\nu }\left( f\right) e^{i\lambda
_{\nu }x}
\end{equation*}%
such that the interval $\left( \frac{\alpha k}{2},\frac{\alpha \left(
k+1\right) }{2}\right) $ does not contain any $\lambda _{\nu }.$ Applying
Lemma 1.10.2 of \cite{BML} we easily verify that 
\begin{equation*}
S_{k}^{\ast }f\left( x\right) -f\left( x\right) =\int_{0}^{\infty }\varphi
_{x}\left( t\right) \Psi _{k}\left( t\right) dt,
\end{equation*}%
where\ $\varphi _{x}\left( t\right) :=f\left( x+t\right) +f\left( x-t\right)
-2f\left( x\right) $ and $\Psi _{k}\left( t\right) =\Psi _{\frac{\alpha k}{2}%
,\frac{\alpha \left( k+1\right) }{2}}\left( t\right) ,$ i.e. 
\begin{equation*}
\Psi _{k}\left( t\right) =\frac{4\sin \frac{\alpha t}{4}\sin \frac{\alpha
\left( 2k+1\right) t}{4}}{\alpha \pi t^{2}}
\end{equation*}%
(see also \cite{ASB}, p.41). Evidently, if the interval $\left( \frac{\alpha
k}{2},\frac{\alpha \left( k+1\right) }{2}\right) $ contains a Fourier
exponent $\lambda _{\nu },$ then 
\begin{equation*}
S_{\frac{\alpha k}{2}}f\left( x\right) =S_{k+1}^{\ast }f\left( x\right)
-\left( A_{\nu }\left( f\right) e^{i\lambda _{\nu }x}+A_{-\nu }\left(
f\right) e^{-i\lambda _{\nu }x}\right) .
\end{equation*}%
Since (see \cite[p.78]{ABI} and \cite[p. 7]{ADB}) 
\begin{equation*}
\left\{ \sum_{\nu =-\infty }^{\infty }\left\vert A_{\nu }\left( f\right)
\right\vert ^{q}\right\} ^{1/q}\leq \left\Vert f\right\Vert _{B^{p}}\text{ \
\ and \ \ }\left\Vert f\right\Vert _{B^{p}}\leq \left\Vert f\right\Vert
_{S^{p}}\text{,}
\end{equation*}%
where $\left\Vert \cdot \right\Vert _{B^{p,}}$\ with $p\geq 1,$ is the
Besicovitch norm,\ so we have 
\begin{equation*}
\left\vert A_{\pm \nu }\left( f\right) \right\vert =\left\vert A_{\pm \nu
}\left( f-g_{\alpha \mu /2}\right) \right\vert \leq \left\Vert f-g_{\alpha
\mu /2}\right\Vert _{S^{p}}=E_{\alpha \mu /2}\left( f\right) _{S^{p}},
\end{equation*}%
for some\ $g_{\alpha \mu /2}\in B_{\alpha \mu /2},$ with $\alpha k/2<\alpha
\mu /2<\lambda _{\nu }.$ Therefore, the deviation 
\begin{equation*}
\left\{ \sum_{k=0}^{n}a_{n,k}\left\vert S_{\frac{\alpha k}{2}}f\left(
x\right) -f\left( x\right) \right\vert ^{q}\right\} ^{1/q}
\end{equation*}%
can be estimated from above by 
\begin{equation*}
\left\{ \sum_{k=0}^{n}a_{n,k}\left\vert \int_{0}^{\infty }\varphi _{x}\left(
t\right) \Psi _{k+\kappa }\left( t\right) dt\right\vert ^{q}\right\}
^{1/q}+\left\{ \sum_{k=0}^{n}a_{n,k}\left( E_{\alpha k/2}\left( f\right)
_{S^{p}}\right) ^{q}\right\} ^{1/q},
\end{equation*}%
where $\kappa $ equals $0$ or $1$. Applying the Minkowski inequality we
obtain 
\begin{equation*}
\left\{ \sum_{k=0}^{n}a_{n,k}\left\vert \int_{0}^{\infty }\varphi _{x}\left(
t\right) \Psi _{k+\kappa }\left( t\right) dt\right\vert ^{q}\right\} ^{1/q}
\end{equation*}%
\begin{equation*}
=\left\{ \sum_{k=0}^{n}a_{n,k}\left\vert \left( \int\limits_{0}^{\frac{2\pi 
}{\alpha }a_{n,n}}+\int\limits_{\frac{2\pi }{\alpha }a_{n,n}}^{\frac{2\pi }{%
\alpha }}+\int\limits_{\frac{2\pi }{\alpha }}^{\infty }\right) \varphi
_{x}\left( t\right) \Psi _{k+\kappa }\left( t\right) dt\right\vert
^{q}\right\} ^{1/q}
\end{equation*}%
\begin{equation*}
\leq \left\{ \sum_{k=0}^{n}a_{n,k}\left\vert I_{1}(k)\right\vert
^{q}\right\} ^{1/q}+\left\{ \sum_{k=0}^{n}a_{n,k}\left\vert
I_{2}(k)\right\vert ^{q}\right\} ^{1/q}+\left\{
\sum_{k=0}^{n}a_{n,k}\left\vert I_{3}(k)\right\vert ^{q}\right\} ^{1/q}.
\end{equation*}%
By (\ref{1}), integrating by parts, we obtain%
\begin{equation*}
\left\{ \sum_{k=0}^{n}a_{n,k}\left\vert I_{1}(k)\right\vert ^{q}\right\}
^{1/q}\leq \left\{ \sum_{k=0}^{n}a_{n,k}\left\vert \frac{4}{\alpha \pi }%
\int\limits_{0}^{\frac{2\pi }{\alpha }a_{n,n}}\varphi _{x}\left( t\right) 
\frac{\sin \frac{\alpha t}{4}}{t^{2}}\sin \frac{\alpha t}{4}(2k+2\kappa
+1)dt\right\vert ^{q}\right\} ^{1/q}
\end{equation*}%
\begin{equation*}
\leq \frac{1}{\pi }\int\limits_{0}^{\frac{2\pi }{\alpha }a_{n,n}}\frac{%
\left\vert \varphi _{x}\left( t\right) \right\vert }{t}dt=\frac{1}{\pi }%
\int\limits_{0}^{\frac{2\pi }{\alpha }a_{n,n}}\frac{1}{t}\left( \frac{d}{dt}%
\int\limits_{0}^{t}\left\vert \varphi _{x}\left( s\right) \right\vert
ds\right) dt
\end{equation*}%
\begin{equation*}
=\frac{1}{\pi }\left[ \frac{1}{t}\int\limits_{0}^{t}\left\vert \varphi
_{x}\left( s\right) \right\vert ds\right] _{t=0}^{t=\frac{2\pi }{\alpha }%
a_{n,n}}+\frac{1}{\pi }\int\limits_{0}^{\frac{2\pi }{\alpha }a_{n,n}}\frac{1%
}{t^{2}}\left( \int\limits_{0}^{t}\left\vert \varphi _{x}\left( s\right)
\right\vert ds\right) dt
\end{equation*}%
\begin{equation*}
=\frac{1}{\pi }w_{x}f\left( \frac{2\pi }{\alpha }a_{n,n}\right) _{1}+\frac{1%
}{\pi }\int\limits_{0}^{\frac{2\pi }{\alpha }a_{n,n}}\frac{1}{t}%
w_{x}f\left( t\right) _{1}dt
\end{equation*}%
\begin{equation*}
\ll w_{x}f\left( a_{n,n}\right) _{1}+\int\limits_{0}^{\frac{2\pi }{\alpha }%
a_{n,n}}\frac{1}{t}w_{x}f\left( t\right) _{1}dt
\end{equation*}%
\begin{eqnarray}
&=&\frac{1}{\pi }w_{x}f\left( \frac{2\pi }{\alpha }a_{n,n}\right) _{1}+\frac{%
1}{\pi }\int\limits_{0}^{\frac{2\pi }{\alpha }a_{n,n}}\frac{1}{t}%
w_{x}f\left( t\right) _{1}dt.  \notag \\
&\ll &w_{x}f\left( a_{n,n}\right) _{1}+\int\limits_{0}^{\frac{2\pi }{\alpha 
}a_{n,n}}\frac{1}{t}w_{x}f\left( t\right) _{1}dt.
\end{eqnarray}%
It is clear that $w_{x}f\left( \delta \right) _{1}/\delta $ is nondecreasing
with respect to $\delta >0$ and $w_{x}f\left( \delta \right) _{1}\leq
w_{x}f\left( \delta \right) _{p}$ for $p\geq 1$. Using these properties we
have%
\begin{equation*}
\left\{ \sum_{k=0}^{n}a_{n,k}\left\vert I_{1}(k)\right\vert ^{q}\right\}
^{1/q}\ll a_{n,n}\int\limits_{a_{n,n}}^{\pi }\frac{w_{x}f\left( t\right)
_{1}}{t^{2}}+\int\limits_{0}^{a_{n,n}}\frac{1}{t}w_{x}f\left( \frac{2\pi }{%
\alpha }t\right) _{1}dt
\end{equation*}%
\begin{equation*}
\ll \left\{ a_{n,n}\int\limits_{a_{n,n}}^{\pi }\frac{\left( w_{x}f\left(
t\right) _{1}\right) ^{p}}{t^{2}}\right\} ^{\frac{1}{p}}+\int%
\limits_{0}^{a_{n,n}}\frac{1}{t}w_{x}f\left( t\right) _{1}dt.
\end{equation*}%
Since $f\in \Omega _{\alpha ,p}\left( w_{x}\right) $ and (\ref{7}) holds,
Lemma 1 and (\ref{6}) give%
\begin{equation*}
\left\{ \sum_{k=0}^{n}a_{n,k}\left\vert I_{1}(k)\right\vert ^{q}\right\}
^{1/q}=O\left( a_{nn}H_{x}\left( a_{nn}\right) \right) .
\end{equation*}%
If (\ref{5}) holds, then%
\begin{equation*}
a_{n,\mu }-a_{n,m}\leq \left\vert a_{n,\mu }-a_{n,m}\right\vert \leq
\sum\limits_{k=\mu }^{m-1}\left\vert a_{n,k}-a_{n,k+1}\right\vert \leq
Ka_{n,m}
\end{equation*}%
for any $m\geq \mu \geq 0.$ Hence we have 
\begin{equation}
a_{n,\mu }\leq \left( K+1\right) a_{n,m}.  \label{16}
\end{equation}%
From this, we get%
\begin{eqnarray*}
&&\left\{ \sum_{k=0}^{n}a_{n,k}\left\vert I_{2}(k)\right\vert ^{q}\right\}
^{1/q} \\
&\leq &\left\{ \left( K+1\right) a_{n,n}\right\} ^{\frac{1}{q}}\left\{
\sum\limits_{k=0}^{n}\left\vert \frac{4}{\alpha \pi }\int\limits_{\frac{%
2\pi }{\alpha }a_{n,n}}^{\frac{2\pi }{\alpha }}\frac{\varphi _{x}\left(
t\right) \sin \frac{\alpha t}{4}}{t^{2}}\sin \frac{\alpha t}{4}(2k+2\kappa
+1)dt\right\vert ^{q}\right\} ^{\frac{1}{q}}
\end{eqnarray*}%
\begin{equation*}
\ll \frac{8}{\alpha ^{2}}\left( a_{nn}\right) ^{\frac{1}{q}}\left\{
\sum\limits_{k=0}^{n}\left\vert \frac{\alpha }{2\pi }\int\limits_{\frac{%
2\pi }{\alpha }a_{n,n}}^{\frac{2\pi }{\alpha }}\frac{\varphi _{x}\left(
t\right) \sin \frac{\alpha t}{4}}{t^{2}}\sin \frac{\alpha t}{4}(2\kappa
+1)\cos \frac{\alpha kt}{2}dt\right\vert ^{q}\right\} ^{\frac{1}{q}}
\end{equation*}%
\begin{equation*}
+\frac{8}{\alpha ^{2}}\left( a_{nn}\right) ^{\frac{1}{q}}\left\{
\sum\limits_{k=0}^{n}\left\vert \frac{\alpha }{2\pi }\int\limits_{\frac{%
2\pi }{\alpha }a_{n,n}}^{\frac{2\pi }{\alpha }}\frac{\varphi _{x}\left(
t\right) \sin \frac{\alpha t}{4}}{t^{2}}\cos \frac{\alpha t}{4}(2\kappa
+1)\sin \frac{\alpha kt}{2}dt\right\vert ^{q}\right\} ^{\frac{1}{q}}.
\end{equation*}%
Using inequality (\ref{13}), we have%
\begin{equation*}
\left\{ \sum_{k=0}^{n}a_{n,k}\left\vert I_{2}(k)\right\vert ^{q}\right\}
^{1/q}\ll \left( a_{nn}\right) ^{\frac{1}{q}}\left\{ \int\limits_{\frac{%
2\pi }{\alpha }a_{n,n}}^{\frac{2\pi }{\alpha }}\frac{\left\vert \varphi
_{x}\left( t\right) \right\vert ^{p}}{t^{1+\frac{p}{q}}}\right\} ^{\frac{1}{p%
}}.
\end{equation*}%
Integrating by parts, we obtain%
\begin{equation*}
\left\{ \sum_{k=0}^{n}a_{n,k}\left\vert I_{2}(k)\right\vert ^{q}\right\}
^{1/q}=\left( a_{nn}\right) ^{\frac{1}{q}}\left\{ \left[ \frac{1}{t^{1+\frac{%
p}{q}}}\int\limits_{0}^{t}\left\vert \varphi _{x}\left( t\right)
\right\vert ^{p}ds\right] _{t=\frac{2\pi }{\alpha }a_{n,n}}^{t=\frac{2\pi }{%
\alpha }}\right.
\end{equation*}%
\begin{equation*}
\left. \left( 1+\frac{p}{q}\right) \int\limits_{\frac{2\pi }{\alpha }%
a_{n,n}}^{\frac{2\pi }{\alpha }}\frac{1}{t^{2+\frac{p}{q}}}\left(
\int\limits_{0}^{t}\left\vert \varphi _{x}\left( t\right) \right\vert
^{p}ds\right) dt\right\} ^{\frac{1}{p}}
\end{equation*}%
\begin{equation*}
=\left( a_{nn}\right) ^{\frac{1}{q}}\left\{ \left[ \frac{1}{t^{\frac{p}{q}}}%
\left( w_{x}f\left( t\right) _{p}\right) ^{p}\right] _{t=\frac{2\pi }{\alpha 
}a_{n,n}}^{t=\frac{2\pi }{\alpha }}+\left( 1+\frac{p}{q}\right)
\int\limits_{\frac{2\pi }{\alpha }a_{n,n}}^{\frac{2\pi }{\alpha }}\frac{1}{%
t^{1+\frac{p}{q}}}\left( w_{x}f\left( t\right) _{p}\right) ^{p}dt\right\} ^{%
\frac{1}{p}}
\end{equation*}%
\begin{equation}
\ll \left( a_{nn}\right) ^{\frac{1}{q}}\left\{ \left( w_{x}f\left( \frac{%
2\pi }{\alpha }\right) _{p}\right) ^{p}+\int\limits_{\frac{2\pi }{\alpha }%
a_{n,n}}^{\frac{2\pi }{\alpha }}\frac{1}{t^{1+\frac{p}{q}}}\left(
w_{x}f\left( t\right) _{p}\right) ^{p}dt\right\} ^{\frac{1}{p}}.  \label{17}
\end{equation}%
Since $f\in \Omega _{\alpha ,p}\left( w_{x}\right) $, (\ref{6}) gives%
\begin{equation*}
\left\{ \sum_{k=0}^{n}a_{n,k}\left\vert I_{2}(k)\right\vert ^{q}\right\}
^{1/q}\ll \left( a_{nn}\right) ^{\frac{1}{q}}\left\{ \left( w_{x}\left( \pi
\right) \right) ^{p}+\int\limits_{a_{nn}}^{\pi }\frac{1}{t^{1+\frac{p}{q}}}%
\left( w_{x}\left( t\right) \right) ^{p}dt\right\} ^{\frac{1}{p}}
\end{equation*}%
\begin{equation*}
\ll \left\{ \left( a_{nn}\right) ^{\frac{p}{q}}\int\limits_{a_{nn}}^{\pi }%
\frac{\left( w_{x}\left( t\right) \right) ^{p}}{t^{1+\frac{p}{q}}}dt\right\}
^{\frac{1}{p}}=O\left( a_{nn}H_{x}\left( a_{nn}\right) \right) .
\end{equation*}%
For the third term we obtain 
\begin{equation*}
\left\{ \sum_{k=0}^{n}a_{n,k}\left\vert I_{3}(k)\right\vert ^{q}\right\}
^{1/q}\leq
\end{equation*}%
\begin{equation*}
\leq \left\{ \sum_{k=0}^{n}a_{n,k}\left\vert \sum_{\mu =1}^{\infty
}\int\limits_{\frac{2\pi }{\alpha }\mu }^{\frac{2\pi }{\alpha }\left( \mu
+1\right) }\left[ \varphi _{x}\left( t\right) -\Phi _{x}f\left( \delta
_{k},t\right) \right] \Psi _{k+\kappa }\left( t\right) dt\right\vert
^{q}\right\} ^{1/q}
\end{equation*}%
\begin{equation*}
+\left\{ \sum_{k=0}^{n}a_{n,k}\left\vert \sum_{\mu =1}^{\infty
}\int\limits_{\frac{2\pi }{\alpha }\mu }^{\frac{2\pi }{\alpha }\left( \mu
+1\right) }\Phi _{x}f\left( \delta _{k},t\right) \Psi _{k+\kappa }\left(
t\right) dt\right\vert ^{q}\right\} ^{1/q}
\end{equation*}%
\begin{equation*}
=\left\{ \sum_{k=0}^{n}a_{n,k}\left\vert I_{31}(k)\right\vert ^{q}\right\}
^{1/q}+\left\{ \sum_{k=0}^{n}a_{n,k}\left\vert I_{32}(k)\right\vert
^{q}\right\} ^{1/q}
\end{equation*}%
and 
\begin{equation*}
\left\vert I_{31}(k)\right\vert \leq \frac{4}{\alpha \pi }\sum_{\mu
=1}^{\infty }\int\limits_{\frac{2\pi }{\alpha }\mu }^{\frac{2\pi }{\alpha }%
\left( \mu +1\right) }\left\vert \varphi _{x}\left( t\right) -\Phi
_{x}f\left( \delta _{k},t\right) \right\vert t^{-2}dt
\end{equation*}%
\begin{equation*}
\leq \frac{4}{\alpha \pi }\sum_{\mu =1}^{\infty }\int\limits_{\frac{2\pi }{%
\alpha }\mu }^{\frac{2\pi }{\alpha }\left( \mu +1\right) }\left[ \frac{1}{%
\delta _{k}t^{2}}\int_{0}^{\delta _{k}}\left\vert \varphi _{x}\left(
t\right) -\varphi _{x}\left( t+u\right) \right\vert du\right] dt
\end{equation*}%
\begin{equation*}
=\frac{4}{\alpha \pi }\frac{1}{\delta _{k}}\int\limits_{0}^{\delta
_{k}}\sum_{\mu =1}^{\infty }\left\{ \int\limits_{\frac{2\pi }{\alpha }\mu
}^{\frac{2\pi }{\alpha }\left( \mu +1\right) }\frac{1}{t^{2}}\left\vert
\varphi _{x}\left( t\right) -\varphi _{x}\left( t+u\right) \right\vert
dt\right\} du
\end{equation*}%
\begin{equation*}
=\frac{4}{\alpha \pi }\frac{1}{\delta _{k}}\int\limits_{0}^{\delta
_{k}}\sum_{\mu =1}^{\infty }\left\{ \left[ \frac{1}{t^{2}}%
\int_{0}^{t}\left\vert \varphi _{x}\left( s\right) -\varphi _{x}\left(
s+u\right) \right\vert ds\right] _{t=\frac{2\pi }{\alpha }\mu }^{t=\frac{%
2\pi }{\alpha }\left( \mu +1\right) }\right.
\end{equation*}%
\begin{equation*}
+\left. 2\int\limits_{\frac{2\pi }{\alpha }\mu }^{\frac{2\pi }{\alpha }%
\left( \mu +1\right) }\left[ \frac{1}{t^{3}}\int_{0}^{t}\left\vert \varphi
_{x}\left( s\right) -\varphi _{x}\left( s+u\right) \right\vert ds\right]
dt\right\} du
\end{equation*}%
\begin{equation*}
\ll \left\vert \frac{1}{\delta _{k}}\int\limits_{0}^{\delta _{k}}\sum_{\mu
=1}^{\infty }\left\{ \frac{1}{[\frac{2\pi }{\alpha }\left( \mu +1\right)
]^{2}}\int\limits_{0}^{\frac{2\pi }{\alpha }\left( \mu +1\right)
}\left\vert \varphi _{x}\left( s\right) -\varphi _{x}\left( s+u\right)
\right\vert ds\right. \right.
\end{equation*}%
\begin{equation*}
\left. -\left. \frac{1}{[\frac{2\pi }{\alpha }\mu ]^{2}}\int\limits_{0}^{%
\frac{2\pi }{\alpha }\mu }\left\vert \varphi _{x}\left( s\right) -\varphi
_{x}\left( s+u\right) \right\vert ds\right\} du\right\vert
\end{equation*}%
\begin{equation*}
+\frac{1}{\delta _{k}}\int\limits_{0}^{\delta _{k}}\sum_{\mu =1}^{\infty
}\left\{ \int\limits_{\frac{2\pi }{\alpha }\mu }^{\frac{2\pi }{\alpha }%
\left( \mu +1\right) }\left[ \frac{1}{t^{3}}\int_{0}^{t}\left\vert \varphi
_{x}\left( s\right) -\varphi _{x}\left( s+u\right) \right\vert ds\right]
dt\right\} du.
\end{equation*}%
Since $f\in \Omega _{\alpha ,p}\left( w_{x}\right) $, thus for any $x$ 
\begin{eqnarray*}
\lim_{\zeta \rightarrow \infty }\frac{1}{\zeta ^{2}}\int\limits_{0}^{\zeta
}\left\vert \varphi _{x}\left( s\right) -\varphi _{x}\left( s+u\right)
\right\vert ds &\leq &\lim_{\zeta \rightarrow \infty }\frac{1}{\zeta }%
w_{x}\left( u\right) \leq \lim_{\zeta \rightarrow \infty }\frac{1}{\zeta }%
w_{x}\left( \delta _{k}\right) \\
&\leq &\lim_{\zeta \rightarrow \infty }\frac{1}{\zeta }w_{x}\left( \pi
\right) =0,
\end{eqnarray*}%
and therefore 
\begin{equation*}
\left\vert I_{31}(k)\right\vert \leq \frac{1}{\delta _{k}}%
\int\limits_{0}^{\delta _{k}}\frac{\alpha }{2\pi }\left[ \frac{\alpha }{%
2\pi }\int_{0}^{2\pi /\alpha }\left\vert \varphi _{x}\left( s\right)
-\varphi _{x}\left( s+u\right) \right\vert ds\right] du
\end{equation*}%
\begin{equation*}
+\frac{1}{\delta _{k}}\int\limits_{0}^{\delta _{k}}w_{x}\left( u\right)
du\sum_{\mu =1}^{\infty }\left\{ \int\limits_{\frac{2\pi }{\alpha }\mu }^{%
\frac{2\pi }{\alpha }\left( \mu +1\right) }\frac{1}{t^{2}}dt\right\}
\end{equation*}%
\begin{equation*}
\ll \frac{1}{\delta _{k}}\int\limits_{0}^{\delta _{k}}w_{x}\left( u\right)
du+w_{x}\left( \delta _{k}\right) \sum_{\mu =1}^{\infty }\frac{1}{\frac{2\pi 
}{\alpha }\mu ^{2}}\ll w_{x}\left( \delta _{k}\right) .
\end{equation*}%
Next, we will estimate the term $\left\vert I_{32}(k)\right\vert .$ So, 
\begin{equation*}
I_{32}(k)=\frac{2}{\alpha \pi }\sum_{\mu =1}^{\infty }\int\limits_{\frac{%
2\pi }{\alpha }\mu }^{\frac{2\pi }{\alpha }\left( \mu +1\right) }\frac{\Phi
_{x}f\left( \delta _{k},t\right) }{t^{2}}\frac{d}{dt}\left( -\frac{\cos 
\frac{\alpha t\left( k+\kappa \right) }{2}}{\frac{\alpha \left( k+\kappa
\right) }{2}}+\frac{\cos \frac{\alpha t\left( k+\kappa +1\right) }{2}}{\frac{%
\alpha \left( k+\kappa +1\right) }{2}}\right) dt
\end{equation*}%
\begin{equation*}
=\frac{2}{\alpha \pi }\sum_{\mu =1}^{\infty }\left[ \frac{\Phi _{x}f\left(
\delta _{k},t\right) }{t^{2}}\left( -\frac{\cos \frac{\alpha t\left(
k+\kappa \right) }{2}}{\frac{\alpha \left( k+\kappa \right) }{2}}+\frac{\cos 
\frac{\alpha t\left( k+\kappa +1\right) }{2}}{\frac{\alpha \left( k+\kappa
+1\right) }{2}}\right) \right] _{t=\frac{2\pi }{\alpha }\mu }^{t=\frac{2\pi 
}{\alpha }\left( \mu +1\right) }
\end{equation*}%
\begin{eqnarray*}
&&+\frac{2}{\alpha \pi }\sum_{\mu =1}^{\infty }\int\limits_{\frac{2\pi }{%
\alpha }\mu }^{\frac{2\pi }{\alpha }\left( \mu +1\right) }\frac{d}{dt}\left( 
\frac{\Phi _{x}f\left( \delta _{k},t\right) }{t^{2}}\right) \left( \frac{%
\cos \frac{\alpha t\left( k+\kappa \right) }{2}}{\frac{\alpha \left(
k+\kappa \right) }{2}}-\frac{\cos \frac{\alpha t\left( k+\kappa +1\right) }{2%
}}{\frac{\alpha \left( k+\kappa +1\right) }{2}}\right) dt \\
&=&I_{321}\left( k\right) +I_{322}\left( k\right)
\end{eqnarray*}%
Since $f\in \Omega _{\alpha ,p}\left( w_{x}\right) $, thus for any $x$
(using (\ref{W}))%
\begin{equation*}
\lim_{\zeta \rightarrow \infty }\left\vert \frac{\Phi _{x}f\left( \delta
_{k},\frac{2\pi }{\alpha }\zeta \right) }{\left[ \frac{2\pi }{\alpha }\zeta %
\right] ^{2}}\left( -\frac{\cos \left[ \pi \zeta (k+\kappa )\right] }{\frac{%
\alpha \left( k+\kappa \right) }{2}}+\frac{\cos \left[ \pi \zeta \left(
k+\kappa +1\right) \right] }{\frac{\alpha \left( k+\kappa +1\right) }{2}}%
\right) \right\vert
\end{equation*}%
\begin{equation*}
\leq \lim_{\zeta \rightarrow \infty }\frac{w_{x}\left( \delta _{k}\right)
+w_{x}\left( \frac{2\pi }{\alpha }\zeta \right) }{2\pi ^{2}\zeta ^{2}k}\ll
\lim_{\zeta \rightarrow \infty }\frac{w_{x}\left( \delta _{k}\right) +\zeta
w_{x}\left( \frac{2\pi }{\alpha }\right) }{\zeta ^{2}k}\ll w_{x}\left( \pi
\right) \lim_{\zeta \rightarrow \infty }\frac{1+\zeta }{\zeta ^{2}}=0,
\end{equation*}%
and therefore%
\begin{equation*}
I_{321}\left( k\right) =\frac{2}{\alpha \pi }\sum_{\mu =1}^{\infty }\left[ 
\frac{\Phi _{x}f\left( \delta _{k},\frac{2\pi }{\alpha }\left( \mu +1\right)
\right) }{\left[ \frac{2\pi }{\alpha }\left( \mu +1\right) \right] ^{2}}%
\left( -\frac{\cos \left[ \pi \left( \mu +1\right) \left( k+\kappa \right) %
\right] }{\frac{\alpha \left( k+\kappa \right) }{2}}\right. \right.
\end{equation*}%
\begin{equation*}
+\left. \frac{\cos \left[ \pi \left( \mu +1\right) \left( k+\kappa +1\right) %
\right] }{\frac{\alpha \left( k+\kappa +1\right) }{2}}\right)
\end{equation*}%
\begin{equation*}
-\left. \frac{\Phi _{x}f\left( \delta _{k},\frac{2\pi }{\alpha }\mu \right) 
}{\left[ \frac{2\pi }{\alpha }\mu \right] ^{2}}\left( -\frac{\cos \left[ \pi
\mu (k+\kappa )\right] }{\frac{\alpha \left( k+\kappa \right) }{2}}+\frac{%
\cos \left[ \pi \mu \left( k+\kappa +1\right) \right] }{\frac{\alpha \left(
k+\kappa +1\right) }{2}}\right) \right]
\end{equation*}%
\begin{equation*}
=-\frac{2}{\alpha \pi }\frac{\Phi _{x}f\left( \delta _{k},2\pi /\alpha
\right) }{\left[ 2\pi /\alpha \right] ^{2}}\left( -\frac{\left( -1\right)
^{\left( k+\kappa \right) }}{\frac{\alpha \left( k+\kappa \right) }{2}}+%
\frac{\left( -1\right) ^{\left( k+\kappa +1\right) }}{\frac{\alpha \left(
k+\kappa +1\right) }{2}}\right)
\end{equation*}%
\begin{equation*}
=-\frac{1}{\pi ^{3}}\Phi _{x}f\left( \delta _{k},2\pi /\alpha \right) \left(
-1\right) ^{\left( k+\kappa +1\right) }\left( \frac{1}{k+\kappa +1}+\frac{1}{%
k+\kappa }\right) .
\end{equation*}%
Using (\ref{W}), we get 
\begin{equation*}
\left\vert I_{321}\left( k\right) \right\vert \ll \frac{1}{\pi ^{3}}\frac{2}{%
k+1}\left\vert \Phi _{x}f\left( \delta _{k},2\pi /\alpha \right) \right\vert
\leq \frac{2}{\pi ^{3}\left( k+1\right) }\left( w_{x}\left( \delta
_{k}\right) +w_{x}\left( 2\pi /\alpha \right) \right) .
\end{equation*}%
Similarly 
\begin{equation*}
I_{322}\left( k\right) =\frac{2}{\alpha \pi }\sum_{\mu =1}^{\infty
}\int\limits_{\frac{2\pi }{\alpha }\mu }^{\frac{2\pi }{\alpha }\left( \mu
+1\right) }\left( \frac{\frac{d}{dt}\Phi _{x}f\left( \delta _{k},t\right) }{%
t^{2}}-\frac{2\Phi _{x}f\left( \delta _{k},t\right) }{t^{3}}\right)
\end{equation*}%
\begin{equation*}
\cdot \left( \frac{\cos \frac{\alpha t\left( k+\kappa \right) }{2}}{\frac{%
\alpha \left( k+\kappa \right) }{2}}-\frac{\cos \frac{\alpha t\left(
k+\kappa +1\right) }{2}}{\frac{\alpha \left( k+\kappa +1\right) }{2}}\right)
dt
\end{equation*}%
and 
\begin{equation*}
\left\vert I_{322}\left( k\right) \right\vert \ll \frac{8}{\alpha ^{2}\left(
k+1\right) \pi }\sum_{\mu =1}^{\infty }\left[ \int\limits_{\frac{2\pi }{%
\alpha }\mu }^{\frac{2\pi }{\alpha }\left( \mu +1\right) }\frac{\left\vert
\varphi _{x}\left( t+\delta _{k}\right) -\varphi _{x}\left( t\right)
\right\vert }{\delta _{k}t^{2}}dt\right.
\end{equation*}%
\begin{equation*}
+\left. 2\int\limits_{\frac{2\pi }{\alpha }\mu }^{\frac{2\pi }{\alpha }%
\left( \mu +1\right) }\frac{\left\vert \Phi _{x}f\left( \delta _{k},t\right)
\right\vert }{t^{3}}dt\right]
\end{equation*}%
\begin{equation*}
\leq \frac{8}{\alpha ^{2}\left( k+1\right) \pi \delta _{k}}\sum_{\mu
=1}^{\infty }\int\limits_{\frac{2\pi }{\alpha }\mu }^{\frac{2\pi }{\alpha }%
\left( \mu +1\right) }\frac{\left\vert \varphi _{x}\left( t+\delta
_{k}\right) -\varphi _{x}\left( t\right) \right\vert }{t^{2}}dt
\end{equation*}%
\begin{equation*}
+\frac{16}{\alpha ^{2}\left( k+1\right) \pi }\sum_{\mu =1}^{\infty
}\int\limits_{\frac{2\pi }{\alpha }\mu }^{\frac{2\pi }{\alpha }\left( \mu
+1\right) }\frac{w_{x}\left( \delta _{k}\right) +w_{x}\left( t\right) }{t^{3}%
}dt
\end{equation*}%
\begin{equation*}
\ll \frac{1}{\left( k+1\right) \delta _{k}}w_{x}\left( \delta _{k}\right) +%
\frac{1}{k+1}\sum_{\mu =1}^{\infty }\left[ \left( w_{x}\left( \delta
_{k}\right) +w_{x}\left( \frac{2\pi \left( \mu +1\right) }{\alpha }\right)
\right) \frac{\alpha ^{2}}{4\pi ^{2}\mu ^{3}}\right]
\end{equation*}%
\begin{equation*}
\ll w_{x}\left( \delta _{k}\right) +\frac{1}{k+1}\left[ w_{x}\left( \delta
_{k}\right) \sum_{\mu =1}^{\infty }\frac{1}{\mu ^{3}}+\sum_{\mu =1}^{\infty }%
\frac{w_{x}\left( \frac{2\pi \left( \mu +1\right) }{\alpha }\right) }{\mu
^{3}}\right]
\end{equation*}%
\begin{equation*}
\ll w_{x}\left( \delta _{k}\right) +\frac{1}{k+1}\left( w_{x}\left( \delta
_{k}\right) +w_{x}\left( \frac{4\pi }{\alpha }\right) \sum_{\mu =1}^{\infty }%
\frac{\mu +1}{\mu ^{3}}\right)
\end{equation*}%
\begin{equation*}
\ll w_{x}\left( \delta _{k}\right) +\frac{1}{k+1}\left( w_{x}\left( \delta
_{k}\right) +w_{x}\left( \frac{4\pi }{\alpha }\right) \right) .
\end{equation*}%
Therefore 
\begin{equation*}
\left\vert I_{3}\left( k\right) \right\vert \ll w_{x}\left( \delta
_{k}\right) +\frac{1}{k+1}\left( w_{x}\left( \delta _{k}\right) +w_{x}\left( 
\frac{2\pi }{\alpha }\right) +w_{x}\left( \frac{4\pi }{\alpha }\right)
\right)
\end{equation*}%
and thus 
\begin{equation*}
\left\{ \sum_{k=0}^{n}a_{n,k}\left\vert I_{3}(k)\right\vert ^{q}\right\}
^{1/q}\ll \left\{ \sum_{k=0}^{n}a_{n,k}\left( w_{x}\left( \frac{\pi }{k+1}%
\right) +\frac{1}{k+1}w_{x}\left( \frac{\pi }{\alpha }\right) \right)
^{q}\right\} ^{1/q}
\end{equation*}%
\begin{equation*}
\ll \left\{ \sum_{k=0}^{n}a_{n,k}\left( w_{x}\left( \frac{\pi }{k+1}\right)
\right) ^{q}\right\} ^{1/q}.
\end{equation*}%
From (\ref{16}) we obtain%
\begin{equation*}
\sum_{k=0}^{n}a_{n,k}\left( w_{x}\left( \frac{\pi }{k+1}\right) \right)
^{q}\leq \sum\limits_{k=0}^{\left[ \frac{1}{\left( K+1\right) a_{n,n}}\right]
-1}a_{nk}\left( w_{x}\left( \frac{\pi }{k+1}\right) \right) ^{q}
\end{equation*}%
\begin{equation*}
+\sum\limits_{k=\left[ \frac{1}{\left( K+1\right) a_{n,n}}\right]
-1}^{n}a_{nk}\left( w_{x}\left( \frac{\pi }{k+1}\right) \right) ^{q}.
\end{equation*}%
Using (\ref{1}), (\ref{16}) and the monotonicity of the function $w_{x}$,
from (\ref{6}) and (\ref{12}), we get%
\begin{equation*}
\sum_{k=0}^{n}a_{n,k}\left( w_{x}\left( \frac{\pi }{k+1}\right) \right)
^{q}\leq \left( K+1\right) a_{n,n}\sum\limits_{k=0}^{\left[ \frac{1}{\left(
K+1\right) a_{n,n}}\right] -1}\left( w_{x}\left( \frac{\pi }{k+1}\right)
\right) ^{q}
\end{equation*}%
\begin{equation*}
+\left( w_{x}\left( \pi \left( K+1\right) a_{n,n}\right) \right)
^{q}\sum\limits_{k=\left[ \frac{1}{\left( K+1\right) a_{n,n}}\right]
-1}^{n}a_{nk}
\end{equation*}%
\begin{equation*}
\ll a_{n,n}\int\limits_{1}^{\frac{1}{\left( K+1\right) a_{n,n}}}\left(
w_{x}\left( \frac{\pi }{t}\right) \right) ^{q}dt+\left( w_{x}\left(
a_{n,n}\right) \right) ^{q}\ll a_{n,n}\int\limits_{a_{n,n}}^{\pi }\frac{%
\left( w_{x}\left( u\right) \right) ^{q}}{u^{2}}du+\left( w_{x}\left(
a_{n,n}\right) \right) ^{q}
\end{equation*}%
\begin{equation*}
\leq a_{n,n}\int\limits_{a_{n,n}}^{\pi }\frac{\left( w_{x}\left( u\right)
\right) ^{q}}{u^{1+p/q+1-p/q}}du+\left( 4w_{x}\left( \frac{a_{n,n}}{2}%
\right) \right) ^{q}
\end{equation*}%
\begin{equation*}
\leq \left( a_{n,n}\right) ^{p/q}\int\limits_{a_{n,n}}^{\pi }\frac{\left(
w_{x}\left( u\right) \right) ^{q}}{u^{1+p/q}}du+\left( 8\int\limits_{\frac{%
a_{n,n}}{2}}^{a_{n,n}}\frac{w_{x}\left( u\right) }{u}du\right) ^{q}
\end{equation*}%
\begin{equation*}
\ll \left( a_{n,n}\right) ^{p/q}\int\limits_{a_{n,n}}^{\pi }\frac{\left(
w_{x}\left( u\right) \right) ^{q}}{u^{1+p/q}}du+\left(
\int\limits_{0}^{a_{n,n}}\frac{w_{x}\left( u\right) }{u}du\right) ^{q}\ll
\left( a_{n,n}H_{x}\left( a_{n,n}\right) \right) ^{q}.
\end{equation*}%
Summing up we obtain that (\ref{8}) is proved and the proof is complete.

\subsection{Proof of Theorem 2}

Under the notation of the before proof we can write 
\begin{equation*}
\left\{ \sum_{k=0}^{n}a_{n,k}\left\vert \int_{0}^{\infty }\varphi _{x}\left(
t\right) \Psi _{k+\kappa }\left( t\right) dt\right\vert ^{q}\right\} ^{1/q}
\end{equation*}%
\begin{equation*}
=\left\{ \sum_{k=0}^{n}a_{n,k}\left\vert \left( \int\limits_{0}^{\frac{2\pi 
}{\alpha }a_{n,0}}+\int\limits_{\frac{2\pi }{\alpha }a_{n,0}}^{\frac{2\pi }{%
\alpha }}+\int\limits_{\frac{2\pi }{\alpha }}^{\infty }\right) \varphi
_{x}\left( t\right) \Psi _{k+\kappa }\left( t\right) dt\right\vert
^{q}\right\} ^{1/q}
\end{equation*}%
\begin{equation*}
\leq \left\{ \sum_{k=0}^{n}a_{n,k}\left\vert J_{1}(k)\right\vert
^{q}\right\} ^{1/q}+\left\{ \sum_{k=0}^{n}a_{n,k}\left\vert
J_{2}(k)\right\vert ^{q}\right\} ^{1/q}+\left\{
\sum_{k=0}^{n}a_{n,k}\left\vert J_{3}(k)\right\vert ^{q}\right\} ^{1/q},
\end{equation*}%
using the Minkowski inequality. Applying the property of the class $RBVS$
instead of the property of $HBVS$ our proof will be similar to the proof of
Theorem 1.

\end{document}